\documentclass[10pt]{article}

\usepackage{amsmath,amsthm,amsfonts,amssymb}

\newtheorem{theorem}{Theorem}[section]
\newtheorem{corollary}[theorem]{Corollary}

\newtheorem{lemma}[theorem]{Lemma}

\newtheorem{remark}[theorem]{Remark}

\def\cH{\mathcal{H}}
\def\cF{\mathcal{F}}

\def\cS{\mathcal{S}}

\def\bC{\mathbb{C}}

\def\bR{\mathbb{R}}

\begin{document}

\title{A note on intermittency for the fractional heat equation}

\author{Raluca M. Balan\footnote{Corresponding author. Department of Mathematics and Statistics, University of Ottawa,
585 King Edward Avenue, Ottawa, ON, K1N 6N5, Canada. E-mail
address: rbalan@uottawa.ca} \footnote{Research supported by a
grant from the Natural Sciences and Engineering Research Council
of Canada.}\and
Daniel Conus\footnote{Lehigh University, Department of Mathematics, 14 East Packer Avenue, Bethlehem, PA, 18109, USA. E-mail address: daniel.conus@lehigh.edu}}

\date{October 30, 2013}
\maketitle

\begin{abstract}
\noindent The goal of the present note is to study intermittency properties for the solution to the fractional heat equation $$\frac{\partial u}{\partial t}(t,x)  =  -(-\Delta)^{\beta/2} u(t,x) + u(t,x)\dot{W}(t,x), \quad t>0,x \in \bR^d$$ with initial condition bounded above and below, where $\beta \in (0,2]$ and the noise $W$ behaves in time like a fractional Brownian motion of index $H>1/2$, and has a spatial covariance given by the Riesz kernel of index $\alpha \in (0,d)$. As a by-product, we obtain that the necessary and sufficient condition for the existence of the solution is $\alpha<\beta$.
\end{abstract}


\vspace{3mm}


\section{Introduction}

In this article we consider the fractional heat equation
\begin{equation}
\left\{\begin{array}{rcl}
\displaystyle \frac{\partial u}{\partial t}(t,x) & = & -(-\Delta)^{\beta/2} u(t,x) + u(t,x)\dot{W}(t,x), \quad t>0,x \in \bR^d \\[2ex]
\displaystyle u(0,x) & = & u_0(x), \quad x \in \bR^d.
\end{array}\right. \label{heat} 
\end{equation}
where $\beta \in (0,2]$, $(-\Delta)^{\beta/2}$ denotes the fractional power of the Laplacian, and $u_0$ is a deterministic function such that
\begin{equation}
\label{initial-bound}
a \leq u_0(x) \leq b \quad \mbox{for all} \quad x \in \bR^d
\end{equation}
for some constants $b \geq a>0$.
We let $W=\{W(\varphi);\varphi \in \cH\}$ be a zero-mean Gaussian process with covariance $$E(W(\varphi)W(\psi))=\langle \varphi,\psi\rangle_{\cH}.$$
Here $\cH$ is a Hilbert space defined as the completion of the space $C_0^{\infty}(\bR_{+} \times \bR^d)$ of infinitely differentiable functions with compact support on $\bR_{+} \times \bR^d$, with respect to the inner product $\langle \cdot, \cdot \rangle_{\cH}$ defined by:
\begin{equation}
\label{def-cov}
\langle \varphi,\psi \rangle_{\cH}= \alpha_H \int_{(\bR_{+} \times \bR^d)^2} \varphi(t,x)\psi(s,y)|t-s|^{2H-2}|x-y|^{-\alpha}dt\, dx\, ds\, dy,
\end{equation}
where $\alpha_H=H(2H-1)$, $H \in (1/2,1)$ and $\alpha \in (0,d)$.
We denote by $\dot{W}$ the formal derivative of $W$. The noise $W$ is spatially homogeneous with spatial covariance given by the Riesz kernel $f(x)=|x|^{-\alpha}$ and behaves in time like a fractional Brownian motion of index $H$. We refer to \cite{B12-IDAQP,B12,BT10} for more details. 

Let $G(t,x)$ be the fundamental solution of $\frac{\partial u}{\partial t}+(-\Delta)^{\beta/2}u=0$ and
$$w(t,x)=\int_{\bR^d}u_0(y)G(t,x-y)dy$$ be the solution of the equation $\frac{\partial u}{\partial t}+(-\Delta)^{\beta/2}u=0$ with initial condition $u(0,x)=u_0(x)$. Note that
\begin{equation}
\label{G-is-density}
\mbox{$G(t,\cdot)$ is the density of $X_t$}
\end{equation}
where $X=(X_t)_{t \geq 0}$ is a symmetric L\'evy process with values in $\bR^d$. If $\beta=2$, then $X$ coincides with a Brownian motion $B=(B_t)_{t \geq 0}$ in $\bR^d$ with variance $2$. If $\beta<2$, then $X$ is a $\beta$-stable L\'evy process given by $X_t=B_{S_t}$, where $(S_t)_{t \geq 0}$ is a ($\beta/2$)-stable subordinator with L\'evy measure
$$\nu(dx)=\frac{\beta/2}{\Gamma(1-\beta/2)}x^{-\beta/2-1}1_{\{x>0\}}dx.$$
Due to \eqref{initial-bound} and \eqref{G-is-density}, it follows that for all $t>0$ and $x \in \bR^d$,
\begin{equation}
\label{bound-w}
a \leq w(t,x) \leq b.
\end{equation}

There is a rich literature dedicated to the case $H=1/2$, when the noise $W$ is white in time. We refer to \cite{dalang,FK13} for some general properties, and to \cite{FK09,CK12,Chen_Dalang} for intermittency properties of the solution to the heat equation with this type of noise. Different methods have to be used for $H>1/2$, since in this case the noise is not a semi-martingale in time.

In the present article, we follow the approach of \cite{hu-nualart09,BT10} for defining the concept of solution. We say that a process $u=\{u(t,x);t \geq 0,x \in \bR^d\}$ defined on a probability space $(\Omega,\cF,P)$ is a {\em mild solution} of \eqref{heat} if it is square-integrable, adapted with respect to the filtration induced by $W$, and satisfies:
$$u(t,x)=w(t,x)+\int_0^t \int_{\bR^d}G(t-s,x-y)u(s,y)W(\delta s,\delta y),$$
where the stochastic integral is interpreted as the divergence operator of $W$ (see (\cite{nualart06}). 
Using Malliavin calculus techniques, it can be shown that the mild solution (if it exists) is unique and has the Wiener chaos decomposition:
\begin{equation}
\label{series}
u(t,x)=\sum_{n \geq 0}I_n(f_n(\cdot,t,x))
\end{equation}
where $I_n$ denotes the multiple Wiener integral (with respect to $W$) of order $n$, and the kernel $f_n(\cdot,t,x)$ is given by:
\begin{eqnarray*}
\lefteqn{f_n(t_1,x_1, \ldots,t_n,x_n,t,x)=} \\
& & G(t-t_n,x-x_n) \ldots G(t_2-t_1,x_2-x_1)w(t_1,x_1)1_{\{0<t_1<\ldots<t_n<t\}}
\end{eqnarray*}
(see page 303 of \cite{hu-nualart09}). By convention, $f_0(t,x)=w(t,x)$ and $I_0$ is the identity map on $\bR$.

The necessary and sufficient condition for the existence of the mild solution is that the series in \eqref{series} converges in $L^2(\Omega)$, i.e.
\begin{equation}
\label{series2}
S(t,x):=\sum_{n \geq 0}\frac{1}{n!}\alpha_n(t,x)<\infty,
\end{equation}
where
$$\alpha_n(t,x)=n!E|I_n(f_n(\cdot,t,x))|^2=(n!)^2 \|\tilde{f}_n(\cdot,t,x)\|_{\cH^{\otimes n}}^2$$
and $\tilde{f}_n(\cdot,t,x)$ is the symmetrization of $f_n(\cdot,t,x)$ in the $n$ variables $(t_1,x_1), \ldots,\linebreak (t_n,x_n)$. If the solution $u$ exists, then $E|u(t,x)|^2=S(t,x)$. We refer to Section 4.1 of \cite{hu-nualart09} and Section 2 of \cite{BT10} for the details. Note that if $u_0(x)=u_0$ for all $x \in \bR^d$, then the law of $u(t,x)$ does not depend on $x$, and hence $\alpha_n(t,x)=\alpha_n(t)$. 

The goal of the present work is to give an upper bound for the $p$-th moment of the solution of \eqref{heat} (for $p \geq 2$), and a lower bound for its second moment. In particular, this will show that, if $u_0(x)$ does not depend on $x$, then the solution $u$ of \eqref{heat} is {\em weakly $\rho$-intermittent}, in a sense which has been recently introduced in \cite{BC}, i.e.
$\gamma_{\rho}(2)>0$ and $\gamma_{\rho}(p)<\infty$ for all $p \geq 2$,
where
$$\gamma_{\rho}(p)=\limsup_{t \to \infty}\frac{1}{t^{\rho}}\log E|u(t,x)|^p$$
is a modified Lyapunov exponent (which does not depend on $x$), and
\begin{equation}
\label{def-rho}
\rho=\frac{2H\beta-\alpha}{\beta-\alpha}.
\end{equation}

As a by-product, we obtain that the necessary and sufficient condition for the existence of the solution is $\alpha<\beta$. Note that this condition is equivalent to
\begin{equation}
\label{DC}
I_{\beta}(\mu):=\int_{\bR^d}\left(\frac{1}{1+|\xi|^{2}}\right)^{\beta/2}\mu(d\xi)
<\infty
\end{equation}
with $\mu(d\xi)=c_{\alpha,d}|\xi|^{-d+\alpha}d\xi$, which is encountered in the study of equations with white noise in time. When $\beta=2$, \eqref{DC} is called {\em Dalang's condition} (see \cite{dalang}).

\section{The result}

The goal of the present article is to prove the following result.

\begin{theorem}
\label{main-th}
The necessary and sufficient condition for equation \eqref{heat} to have a mild solution is $\alpha<\beta$. If the solution $u=\{u(t,x);t \geq 0,x \in \bR^d\}$ exists, then for any $p \geq 2$, for any $x \in \bR^d$ and for any $t>0$ such that
$p t^{2H-\alpha/\beta}>t_1$
$$E|u(t,x)|^p \leq b^p \exp(C_1 p^{(2\beta-\alpha)/(\beta-\alpha)} t^{\rho})$$
and for any $x \in \bR^d$ and for any $t>t_2$,
$$E|u(t,x)|^2 \geq a^2 \exp(C_2 t^{\rho}),$$
where $\rho$ is given by \eqref{def-rho}, $a,b$ are the constants given by \eqref{initial-bound}, and $t_1,t_2,C_1,C_2$ are some positive constants depending on $d,\alpha,\beta$ and $H$.
\end{theorem}

Before giving the proof, we recall from \cite{BT10} that
\begin{equation}
\label{def-alpha-n}
\alpha_n(t,x)=\alpha_H^n \int_{[0,t]^{2n}} \prod_{j=1}^{n}|t_j-s_j|^{2H-2}\psi_{n}({\bf t},{\bf s})d{\bf t}d{\bf s}
\end{equation}
where
$$\psi_{n}({\bf t},{\bf s})=\int_{\bR^{2nd}}\prod_{j=1}^{n}|x_j-y_j|^{-\alpha}
\tilde{f}_n(t_1,x_1,\ldots,t_n,x_n,t,x)
\tilde{f}_n(s_1,y_1, \ldots,s_n,y_n,t,x)d{\bf x}d{\bf y}$$
and we denote ${\bf t}=(t_1, \ldots,t_n)$, ${\bf s}=(s_1, \ldots,s_n)$ with $t_i,s_i \in [0,t]$ and ${\bf x}=(x_1, \ldots,x_n)$, ${\bf y}=(y_1, \ldots,y_n)$ with $x_i,y_i \in \bR^d$.

Note that the Fourier transform of $G(t,\cdot)$ is given by:
\begin{equation}
\label{Fourier-G}\cF G(t,\cdot)(\xi):=\int_{\bR^d}e^{-i \xi \cdot x}G(t,x)dx=\exp(-t|\xi|^{\beta}), \quad \xi \in \bR^d
\end{equation}
where $|\cdot|$ denotes the Euclidean norm in $\bR^d$. Recall that
for any $\varphi, \psi \in L^1(\bR^d)$,
\begin{equation}
\label{harmonic1}
\int_{\bR^d}\int_{\bR^d}\varphi(x)\psi(y)|x-y|^{-\alpha}dxdy=c_{\alpha,d}\int_{\bR^d}\cF \varphi(\xi)\overline{\cF\psi(\xi)}|\xi|^{-d+\alpha}d\xi
\end{equation}
where $\cF \varphi$ is the Fourier transform of $\varphi$, $c_{\alpha,d}=(2\pi)^{-d}C_{\alpha,d}$ and $C_{\alpha,d}$ is the constant given by \eqref{Riesz-const}
(see Appendix A). This identity can be extended to functions $\varphi,\psi \in L^1(\bR^{nd})$:
\begin{eqnarray}
\label{harmonic}
\lefteqn{\int_{\bR^{nd}}\int_{\bR^{nd}}\varphi({\bf x})\psi({\bf y})\prod_{j=1}^{n}|x_j-y_j|^{-\alpha}d{\bf x}d{\bf y}=}\\
\nonumber
& & c_{\alpha,d}^n\int_{\bR^{nd}}\cF \varphi(\xi_1, \ldots,\xi_n)\overline{\cF\psi(\xi_1, \ldots,\xi_n)}\prod_{j=1}^{n}|\xi_j|^{-d+\alpha}d\xi_1 \ldots \xi_n.
\end{eqnarray}

We will use the following elementary inequality.

\begin{lemma}
\label{elem-ineq}
For any $t>0$ and $\eta \in \bR^d$
$$\int_{\bR^d}e^{-t|\xi|^{\beta}}|\xi-\eta|^{-d+\alpha}d\xi \leq K_{d,\alpha,\beta}t^{-\alpha/\beta}$$
where
$$K_{d,\alpha,\beta}:=\sup_{\eta \in \bR^d}\int_{\bR^d}\frac{1}{1+|\xi-\eta|^{\beta}}|\xi|^{-d+\alpha}d\xi.$$
\end{lemma}

\noindent {\bf Proof:} Using the change of variable $z=t^{1/\beta}(\eta-\xi)$, we have:
$$\int_{\bR^d}e^{-t|\xi|^{\beta}}|\xi-\eta|^{-d+\alpha}d\xi = t^{-\alpha/\beta}\int_{\bR^d}e^{-|z-t^{1/\beta}\eta|^{\beta}}|z|^{-d+\alpha}dz.$$
The result follows using the inequality $e^{-x} \leq 1/(1+x)$ for $x>0$.
$\Box$

\vspace{3mm}

\noindent {\bf Proof of Theorem \ref{main-th}:}
{\em Step 1. (Sufficiency and upper bound for the second moment)} Suppose that $\alpha<\beta$. We will prove that the series \eqref{series2} converges, by providing upper bounds for $\psi_{n}({\bf t},{\bf s})$ and $\alpha_n(t,x)$.

By the Cauchy-Schwarz inequality, $\psi_n({\bf t},{\bf s}) \leq \psi_n({\bf t},{\bf t})^{1/2}\psi_n({\bf s},{\bf s})^{1/2}$. So it is enough to consider
the case ${\bf t}={\bf s}$. Let $u_j=t_{\rho(j+1)}-t_{\rho(j)}$ where $\rho$ is a permutation of $\{1, \ldots,n\}$ such that $t_{\rho(1)}<\ldots<t_{\rho(n)}$ and $t_{\rho(n+1)}=t$. Using \eqref{bound-w}, \eqref{Fourier-G} and \eqref{harmonic}, and arguing as in the proof of Lemma 3.2 of \cite{B12}, we obtain:
$$\psi_n({\bf t},{\bf t}) \leq b^2 c_{\alpha,d}^n\int_{\bR^{d}} d\eta_1 \exp(-u_1|\eta_1|^{\beta}) |\eta_1|^{-d+\alpha} \int_{\bR^d}d\eta_2\exp(-u_2|\eta_2|^{\beta})|\eta_2-\eta_1|^{-d+\alpha}$$
 $$\ldots \int_{\bR^d}d\eta_n\exp(-u_n|\eta_n|^{\beta})|\eta_n-\eta_{n-1}|^{-d+\alpha}.$$

\noindent By Lemma \ref{elem-ineq}, it follows that:
$$\psi_n({\bf t},{\bf t}) \leq b^2 c_{\alpha,d}^n K_{d,\alpha,\beta}^n(u_1 \ldots u_n)^{-\alpha/\beta}.$$
By inequality \eqref{ineq-I-beta} (Appendix A),
$K_{d,\alpha,\beta}\leq c_{\beta} I_{d,\alpha,\beta}$,
where $c_{\beta}=2^{\beta/2-1}$ and
$$I_{d,\alpha,\beta}:=\int_{\bR^d}
\left(\frac{1}{1+|\xi|^2}\right)^{\beta/2}|\xi|^{-d+\alpha}d\xi=
\frac{(2\pi)^d c_d\Gamma((\beta-\alpha)/2)\Gamma(\alpha/2)}{2\Gamma(\beta/2)}$$
(see relation \eqref{identity1} and Remark \ref{calculation-I}, Appendix A). Hence,
$$\psi_n({\bf t},{\bf s}) \leq b^2 C_{d,\alpha,\beta}^n [\beta({\bf t}) \beta({\bf s})]^{-\alpha/(2\beta)}$$
where $\beta({\bf t})=u_1 \ldots u_n$, $\beta({\bf s})$ is defined similarly, and $C_{d,\alpha,\beta}>0$ is a constant depending on $d,\alpha,\beta$.
Similarly to the proof of Proposition 3.5 of \cite{BT10}, we have:
\begin{equation}
\label{estimate-alpha}
\alpha_n(t,x) \leq b^2 C_{d,\alpha,\beta,H}^{n} (n!)^{\alpha/\beta} t^{n(2H-\alpha/\beta)},
\end{equation}
where $C_{d,\alpha,\beta,H}>0$ is a constant depending on $d,\alpha,\beta,H$.
Since $\alpha<\beta$, it follows that the series \eqref{series2} converges and
$$E|u(t,x)|^2=\sum_{n \geq 0}\frac{1}{n!}\alpha_n(t,x) \leq b^2 \sum_{n \geq 0}\frac{C_{d,\alpha,\beta,H}^{n}}{(n!)^{1-\alpha/\beta}}t^{n(2H-\alpha/\beta)}\leq b^2 \exp(C_0 t^{\rho}),$$
for all $t>t_0$, where $C_0>0$ and $t_0>0$ are constants depending in $d,\alpha,\beta,H$.
We used the fact that for any $a>0$ and $x>0$,
\begin{equation}
\label{Mittag-Leffler}
\sum_{n \geq 0}\frac{x^n}{(n!)^a} \leq \exp(c_0 x^{1/a}) \quad \mbox{for all} \quad x>x_0,
\end{equation}
where $x_0>0$ and $c_0>0$ are some constants depending on $a$.

\vspace{3mm}

{\em Step 2. (Upper bound for the $p$-the moment)} Note that $u(t,x)=\sum_{n \geq 0}J_n(t,x)$ in $L^2(\Omega)$, where $J_n(t,x)$ lies in the $n$-th order Wiener chaos $\cH_n$ associated to the Gaussian process $W$ (see \cite{nualart06}). Hence,
$$E|u(t,x)|^2=\sum_{n \geq 0}E|J_n(t,x)|^2=\sum_{n \geq 0}\frac{1}{n!}\alpha_n(t,x).$$
We denote by $\|\cdot\|_p$ the $L^p(\Omega)$-norm. We use the fact that for a {\em fixed} Wiener chaos $\cH_n$, the $\|\cdot\|_p$ are equivalent, for all $p \geq 2$ (see the last line of page 62 of \cite{nualart06} with $q=p$ and $p=2$). Hence,
\begin{eqnarray*}
\|J_n(t,x)\|_p & \leq & (p-1)^{n/2}\|J_n(t,x)\|_{2}=(p-1)^{n/2} \left( \frac{1}{n!}\alpha_n(t,x)\right)^{1/2} \\
& \leq & b [(p-1)C_{d,\alpha,\beta,H}]^{n/2} \frac{1}{(n!)^{(\beta-\alpha)/(2\beta)}}
t^{n(2H\beta-\alpha)/(2\beta)}
\end{eqnarray*}
using \eqref{estimate-alpha} for the last inequality. Using Minkowski's inequality for integrals (see Appendix A.1 of \cite{stein70}) and inequality \eqref{Mittag-Leffler}, we obtain that:
$$\|u(t,x)\|_{p} \leq \sum_{n \geq 0}\|J_n(t,x)\|_{p} \leq b \exp(C_1 (p-1)^{\beta/(\beta-\alpha)} t^{\rho})$$
if $p t^{2H-\alpha/\beta}>t_1$, where the constants $C_1>0$ and $t_1>0$ depend on $d,\alpha,\beta,H$.

\vspace{3mm}

{\em Step 3. (Necessity and lower bound for the second moment)}
Suppose that equation \eqref{heat} has a mild solution $u$, i.e. the series \eqref{series2} converges. In particular,
\begin{eqnarray*}
\infty>\alpha_1(t,x) & \geq & a^2 \alpha_H \int_{[0,t]^2} \int_{\bR^{2d}}|r-s|^{2H-s}|y-z|^{-\alpha}G(s,y)G(r,z)dydz drds \\
&=& a^2 \alpha_H c_{\alpha,d} \int_{\bR^d}\left(\int_0^t \int_0^t |r-s|^{2H-2}e^{-(r+s)|\xi|^{\beta}}drds \right)|\xi|^{-d+\alpha}d\xi \\
& \geq & a^2 \alpha_H c_{\alpha,d} c_H \int_{\bR^d}\left(\frac{1}{1/t+|\xi|^{\beta}}\right)^{2H}|\xi|^{-d+\alpha}d\xi,
\end{eqnarray*}
where we used \eqref{harmonic1} for the equality and Theorem 3.1 of \cite{B12-IDAQP} for the last inequality. From here, we infer that
\begin{equation}
\label{ineq-alpha-beta-H}
\alpha<2H\beta.
\end{equation}
In particular, this implies that $\alpha<2\beta$.

Note that one can replace $\psi_{n}({\bf t},{\bf s})$ by $\psi_{n}(t{\bf e}-{\bf t},t{\bf e}-{\bf s})$ in the definition \eqref{def-alpha-n} of $\alpha_n(t,x)$, where ${\bf e}=(1, \ldots,1) \in \bR^n$. By Lemma 2.2 of \cite{B09}, we have:
$$\psi_{n}(t{\bf e}-{\bf t},t{\bf e}-{\bf s})=
E\left[w(t-t^*,x+X_{t^*}^1)w(t-s^*,x+X_{s^*}^2)
\prod_{j=1}^{n}|X_{t_j}^1-X_{s_j}^{2}|^{-\alpha}\right],$$
where $t^*=\max\{t_1, \ldots,t_n\}$, $s^*=\max\{s_1,\ldots,s_n\}$ and $X^1,X^2$ are two independent copies of the L\'evy process $X=(X_t)_{t \geq 0}$ mentioned in the Introduction. (Lemma 2.2 of \cite{B09} was proved for $\beta=2$. The same proof is valid for $\beta<2$.)

Due to \eqref{bound-w}, it follows that
\begin{equation}
\label{ineq-Mn}
a^2 M_n(t) \leq \alpha_n(t,x) \leq b^2 M_n(t)
\end{equation}
where
$$M_n(t):= E\left[ \alpha_H^n \int_{[0,t]^{2n}} \prod_{j=1}^{n}|t_j-s_j|^{2H-2}
\prod_{j=1}^{n}|X_{t_j}^1-X_{s_j}^{2}|^{-\alpha}d{\bf t}d{\bf s}\right]=E(L(t)^n)$$
and $L(t)$ is a random variable defined by:
$$L(t):=\alpha_H \int_0^t \int_0^t |r-s|^{2H-2}|X_{r}^{1}-X_{s}^2|^{-\alpha}drds.$$

To prove that $L(t)$ is finite a.s., we show that its mean is finite. Note that $X_r^1-X_s^2 \stackrel{d}{=}X_{r+s}\stackrel{d}{=}(r+s)^{1/\beta}X_1$, and hence
$$E[L(t)]=\alpha_H C_{d,\alpha,\beta} \int_0^t \int_0^t |r-s|^{2H-2}(r+s)^{-\alpha/\beta} drds,$$
where
$$C_{d,\alpha,\beta}:=E|X_1|^{-\alpha}=\frac{c_d C_{\alpha,d}}{\beta }\Gamma
(\alpha/\beta).$$
(The negative moment of the $\beta$-stable random variable $X_1$ can be computed similarly to \eqref{negative-momentZ}, Appendix A.) Due to \eqref{ineq-alpha-beta-H}, it follows that $E[L(t)]<\infty$.

 By \eqref{ineq-Mn}, we have:
\begin{equation}
\label{ineq-exp-L}
a^2 E(e^{L(t)})\leq E|u(t,x)|^2=\sum_{n \geq 0}\frac{1}{n!}\alpha_n(t,x) \leq b^2 E(e^{L(t)}).
\end{equation}

We consider also the random variable
$$\zeta(t):=\int_0^t \int_0^t |X_r^1-X_s^2|^{-\alpha}drds.$$ Since $|r-s|^{2H-2} \geq (2t)^{2H-2}$ for any $r,s \in [0,t]$,
$L(t) \geq \beta_H t^{2H-2}\zeta(t)$, where $\beta_{H}=\alpha_H 2^{2H-2}$. Hence $\zeta(t)$ is finite a.s.

By the self-similarity (of index $1/\beta$) of the processes $X^1$ and $X^2$, it follows that for any $t>0$ and $c>0$,
$$\zeta(t) \stackrel{d}{=} c^{(2\beta-\alpha)/\beta} \zeta(t/c).$$
In particular, for $c=t^{-(2H-2)\beta/(2\beta-\alpha)}$, we obtain that
$$t^{2H-2}\zeta(t) \stackrel{d}{=} \zeta(t^{\delta}), \quad \mbox{with} \quad \delta=\frac{2H\beta-\alpha}{2\beta-\alpha}$$
and for $c=t$, we obtain that $\zeta(t)\stackrel{d}{=}t^{(2\beta-\alpha)/\beta}\zeta(1)$. Hence,
\begin{equation}
\label{relation-L-zeta}
E(e^{L(t)}) \geq E(e^{\beta_H t^{2H-2} \zeta(t)})=E(e^{\beta_H\zeta(t^{\delta})}).
\end{equation}

The asymptotic behavior of the moments of $\zeta(t)$ was investigated in \cite{bass-chen-rosen09}, under the condition $\alpha<2\beta$. More precisely, under this condition, by relation (2.3) of \cite{bass-chen-rosen09}, we know that:
$$\lim_{n \to \infty}\frac{1}{n}\log \left\{\frac{1}{(n!)^{\alpha/\beta}} E[\zeta(1)^n]\right\}=\log\left(\frac{2\beta}{2\beta-\alpha} \right)^{(2\beta-\alpha)/\beta}+\log \gamma,$$
where $\gamma>0$ is a constant depending on $d,\alpha,\beta$. Hence, there exists some $n_1 \geq 1$ such that for all $n \geq n_1$,
$E[\zeta(1)^n] \geq c^{n}(n!)^{\alpha/\beta}$, where $c>0$ is a constant depending on $d,\alpha,\beta$. Consequently, for any $t>0$,
$$E[\zeta(t)^n] \geq  c^n  t^{n(2\beta-\alpha)/\beta}(n!)^{\alpha/\beta} \quad \mbox{for all} \ n \geq n_1.$$
Hence, for any $\theta>0$,
\begin{equation}
\label{Laplace-zeta}
E(e^{\theta \zeta(t)})=\sum_{n \geq 0}\frac{1}{n!}\theta^n E[\zeta(t)^n] \geq \sum_{n \geq n_1}\frac{1}{(n!)^{1-\alpha/\beta}}\theta^n c^n t^{n(2\beta-\alpha)/\beta}.
\end{equation}
Using \eqref{ineq-exp-L}, \eqref{relation-L-zeta} and \eqref{Laplace-zeta}, we obtain that:
$$\infty>E|u(t,x)|^2 \geq a^2 E(e^{L(t)}) \geq a^2E\left(e^{\beta_H \zeta(t^{\delta})}\right)\geq a^2 \sum_{n \geq n_1} \frac{\beta_H^n c^n t^{n(2H\beta-\alpha)/\beta}}{(n!)^{1-\alpha/\beta}}.$$
{\em This implies that $\alpha<\beta$.} For any $x>0$ and $h \in (0,1)$,
we note that $$E_h(x):=\sum_{n \geq 0}\frac{x^n}{(n!)^h} \geq \left(\sum_{n \geq 0} \frac{(x^{1/h})^n}{n!} \right)^h=\exp(hx^{1/h}).$$
We denote $x_t=\theta c t^{(2\beta-\alpha)/\beta}$ and $h=1-\alpha/\beta$. Writing the last sum in \eqref{Laplace-zeta} as the sum for all terms $n \geq 0$, minus the sum $S_t$ with terms $n\leq n_1$, we see that for all $\theta>0$, and for all $t \geq t_0$,
\begin{eqnarray*}
E(e^{\theta \zeta(t)}) & \geq &
E_{h}(x_t)-S_t \geq \exp(h x_t^{1/h})-S_t \geq \frac{1}{2} \exp(h x_{t}^{1/h}) \\
& \geq & \exp(c_0 \theta^{\beta/(\beta-\alpha)}t^{(2\beta-\alpha)/(\beta-\alpha)}),
\end{eqnarray*}
where $c_0=hc^{1/h}$ and $t_0>0$ is a constant depending on $\theta,\alpha,\beta$. Using this last inequality with $\theta=\beta_H$ and $t^{\delta}$ instead of $t$, we obtain that:
$$E|u(t,x)|^2 \geq a^2 E\left(e^{\beta_{H}\zeta(t^{\delta})}\right) \geq a^2\exp(C_2 t^{\rho}),$$
where $C_2=c_0 \beta_H^{\beta/(\beta-\alpha)}$ depends on $d,\alpha,\beta,H$.
$\Box$

\appendix

\section{Some useful identities}

In this section, we give a result which was used in the proof of Theorem \ref{main-th} for finding an upper bound for $\psi_n({\bf t},{\bf t})$.  This result may be known, but we were not able to find a reference. We state it in a general context.

Following Definition 5.1 of \cite{khoshnevisan-xiao09}, we say that a function $f:\bR^d \to [0,\infty]$ is a kernel of {\em positive type} if it is locally integrable and its Fourier transform in $\cS'(\bR^d)$ is a function $g$ which is non-negative almost everywhere. Here we denote by $\cS'(\bR^d)$ the dual of the space $\cS(\bR^d)$ of rapidly decreasing, infinitely differentiable functions on $\bR^d$.

The Riesz kernel defined by $f(x)=|x|^{-\alpha}$ for $x \in \bR^d \verb2\2 \{0\}$ and $f(0)=\infty$ (with $\alpha \in (0,d)$), is a kernel of positive type. Its Fourier transform in $\cS'(\bR^d)$ is given by $g(\xi)=C_{\alpha,d}|\xi|^{-(d-\alpha)}$ where
\begin{equation}
\label{Riesz-const}C_{\alpha,d}=\pi^{-d/2}2^{-\alpha} \frac{\Gamma((d-\alpha)/2)}{\Gamma(\alpha/2)}
\end{equation}
(see Lemma 1, page 117 of \cite{stein70}).

 Let $f$ be a continuous symmetric kernel of positive type such that $f(x)<\infty$ if and only if $x \not=0$. By Lemma 5.6 of \cite{khoshnevisan-xiao09}, for any Borel probability measures $\mu$ and $\nu$ on $\bR^d$, we have:
 $$\int_{\bR^d} \int_{\bR^d}f(x-y)\mu(dx)\nu(dy)=\frac{1}{(2\pi)^d} \int_{\bR^d}\cF \mu(\xi) \overline{\cF \nu(\xi)}g(\xi)d\xi,$$
where $\cF \mu,\cF \nu$ denote the Fourier transforms of $\mu,\nu$. In particular, if $\mu(dx)=\varphi(x)dx$ and $\nu(dy)=\psi(y)dy$ for some density functions $\varphi,\psi$ in $\bR^d$, then
\begin{equation}
\label{energy-identity}
\int_{\bR^d} \int_{\bR^d}f(x-y)\varphi(x)\psi(y)dxdy=\frac{1}{(2\pi)^d} \int_{\bR^d}\cF \varphi(\xi) \overline{\cF \psi(\xi)}g(\xi)d\xi.
\end{equation}
This relation holds for arbitrary non-negative functions $\varphi,\psi \in L^1(\bR^d)$. (To see this, we consider the normalized functions $\varphi/\|\varphi\|_1$ and  $\psi/\|\psi\|_1$, where $\|\cdot\|_1$ denotes the $L^1(\bR^d)$-norm.)
Using the decomposition $\varphi=\varphi^{+}-\varphi^{-}$ with non-negative functions $\varphi^{+},\varphi^{-}$, we see that \eqref{energy-identity} holds for any functions $\varphi,\psi \in L^1(\bR^d)$. In fact, \eqref{energy-identity} holds for any functions $\varphi,\psi \in L_{\bC}^1(\bR^d)$, replacing $\psi(y)$ by its conjugate $\overline{\psi(y)}$ on the left-hand side. (To see this, we write $\varphi={\varphi}_1+i\varphi_2$ where $\varphi_1,\varphi_2$ are the real and imaginary parts of $\varphi$.)

We consider the Bessel kernel (in $\bR^d$) of order $\beta>0$:
$$G_{d,\beta}(x)=\frac{1}{\Gamma(\beta/2)}\int_{0}^{\infty}u^{\beta/2-1}e^{-u}
\frac{1}{(4\pi u)^{d/2}}e^{-|x|^2/(4u)}du.$$
 Note that $G_{d,\beta}$ is a density function (see Remark \ref{calculation-I} below) and
\begin{equation}
\label{Fourier-Bessel}
\cF G_{d,\beta}(\xi)=\left(\frac{1}{1+|\xi|^2} \right)^{\beta/2}, \quad \xi \in \bR^d.
\end{equation}
Moreover, $G_{d,\alpha}*G_{d,\beta}=G_{d,\alpha+\beta}$ for any $\alpha,\beta>0$
 (see pages 130-135 of \cite{stein70}).

The following result is an extension of relations (3.4) and (3.5) of \cite{dalang-mueller03} to the case of arbitrary $\beta>0$.

\begin{lemma}
\label{lemmaA}
Let $f$ be a continuous symmetric kernel of positive type such that $f(x)<\infty$ if and only if $x\not=0$. Let $\mu(d\xi)=(2\pi)^{-d}g(\xi)d\xi$, where $g$ is the Fourier transform of $f$ in $\cS'(\bR^d)$. Let $\beta>0$ be arbitrary. Then
\begin{equation}
\label{identity1}
\int_{\bR^d}G_{d,\beta}(x)f(x)dx=\int_{\bR^d}
\left(\frac{1}{1+|\xi|^2}\right)^{\beta/2}\mu(d\xi):=I_{\beta}(\mu).
\end{equation}

If $I_{\beta}(\mu)<\infty$, then, for any $a\in \bR^d$,
\begin{equation}
\label{identity2}
\int_{\bR^d}e^{ia \cdot x}G_{d,\beta}(x)f(x)dx=\int_{\bR^d}\left(\frac{1}{1+|\xi-a|^2} \right)^{\beta/2}\mu(d\xi).
\end{equation}
\end{lemma}

\noindent {\bf Proof:} Relation \eqref{identity1} follows from \eqref{energy-identity} with $\varphi=\psi=G_{d,\beta/2}$. On the left-hand side (LHS), we use the fact that $G_{d,\beta/2}*G_{d,\beta/2}=G_{d,\beta}$. On the right-hand side (RHS), we use \eqref{Fourier-Bessel} (with $\beta/2$ instead of $\beta$).

To prove \eqref{identity2}, we apply \eqref{energy-identity} to the complex-valued functions: $$\varphi(x)=\psi(x)=e^{ia \cdot x}G_{d,\beta/2}(x).$$
The term on the LHS is
$$\int_{\bR^d}\int_{\bR^d} e^{ia\cdot (x-y)} G_{d,\beta/2}(x)G_{d,\beta/2}(y)f(x-y)dxdy=
\int_{\bR^d}e^{ia \cdot x}f(x) G_{d,\beta}(x)dx,$$
using Fubini's theorem. The application of Fubini's theorem is justified since
$$\int_{\bR^d}\int_{\bR^d}|e^{ia\cdot (x-y)} G_{d,\beta/2}(x)G_{d,\beta/2}(y)f(x-y)|dxdy=\int_{\bR^d}G_{d,\beta}(x)f(x)dx<\infty.$$
For the term on the RHS, we use the fact that
$$\cF\varphi(\xi)=\int_{\bR^d}e^{-i (\xi-a)\cdot x} G_{d,\beta/2}(x)dx=\cF G_{d,\beta/2}(\xi-a)=\left(\frac{1}{1+|\xi-a|^2}\right)^{\beta/4}.$$
$\Box$

\begin{corollary}
\label{corol-lemmaA}
Let $(f,\mu)$ be as in Lemma \ref{lemmaA} and $\beta>0$ be arbitrary.
Assume that
$I_{\beta}(\mu)<\infty$. Then $$\sup_{a \in \bR^d} \int_{\bR^d}\left(\frac{1}{1+|\xi-a|^2}\right)^{\beta/2}\mu(d\xi)=I_{\beta}(\mu).$$
Consequently,
\begin{equation}
\label{ineq-I-beta}
\sup_{a \in \bR^d} \int_{\bR^d}\frac{1}{1+|\xi-a|^{\beta}}\mu(d\xi) \leq c_{\beta}I_{\beta}(\mu),
\end{equation}
where $c_{\beta}=2^{\beta/2-1}$.
\end{corollary}

\noindent {\bf Proof:} The fact that $I_{\beta}(\mu)$ is smaller than the supremum is obvious. To prove the other inequality, we take absolute values on both sides of \eqref{identity2} and we use the fact that $|\int \cdots |\leq \int |\cdots|$. For the last statement, we use the fact that
$(1+|\xi-a|^2)^{\beta/2} \leq c_{\beta}(1+|\xi-a|^{\beta})$.
$\Box$

\begin{remark}
\label{calculation-I}
{\rm The Bessel kernel $G_{d,\beta}(x)$ arises in statistics as the density of the random vector $X$ given by the following hierarchical model:
$$X|U=u \sim N_{d}(0,2uI) \quad \quad \quad U \sim {\rm Gamma}(\beta/2,1)$$
where $N_d(0,2uI)$ denotes the $d$-dimensional normal distribution with covariance matrix $2uI$, $I$ being the identity matrix. Hence, the term on the LHS of \eqref{identity1} is
$$\int_{\bR^d}G_{d,\beta}(x)f(x)dx=E[f(X)]=
\frac{1}{\Gamma(\beta/2)}\int_{0}^{\infty}u^{\beta/2-1}e^{-u}E[f(X)|U=u]du.$$

This can be computed explicitly if $f(x)=|x|^{-\alpha}$ with $\alpha \in (0,d)$. First, note that if $Z \sim N_{d}(0,2t I)$, then its negative moment of order $-\alpha$ is:
\begin{equation}
\label{negative-momentZ}
E(|Z|^{-\alpha})=\frac{1}{2} C_{\alpha,d}c_d\Gamma(\alpha/2)t^{-\alpha/2}
\end{equation}
where $c_d=2\pi^{d/2}/\Gamma(d/2)$ is the surface area of the unit sphere in $\bR^d$. To see this, we use the fact that $\cF f(\xi)=C_{\alpha,d}|\xi|^{-d+\alpha}d\xi$ in $\cS'(\bR^d)$. Hence,
\begin{eqnarray*}
E(|Z|^{-\alpha})&=& \int_{\bR^d} |x|^{-\alpha}\frac{1}{(4\pi t)^{d/2}}e^{-|x|^2/(4t)}dx=C_{\alpha,d}\int_{\bR^d} |\xi|^{-d+\alpha}e^{-t|\xi|^2}d\xi
\end{eqnarray*}
and \eqref{negative-momentZ} follows by passing to the polar coordinates.
We obtain that
$$\int_{\bR^d}G_{d,\beta}|x|^{-\alpha}dx=
\frac{c_{\alpha,d}c_d\Gamma(\alpha/2)}{2\Gamma(\beta/2)}
\int_0^{\infty} u^{(\beta-\alpha)/2-1}e^{-u}du=
\frac{C_{\alpha,d}c_d\Gamma((\beta-\alpha)/2)\Gamma(\alpha/2)}{2\Gamma(\beta/2)}.$$
(Note that the integral is finite if and only if $\alpha<\beta$.) }
\end{remark}

\end{document}